\documentstyle{amsart}  %

\newcommand{\mathcal}[1]{\cal{#1}}
\newcommand{\mathbb}[1]{\Bbb{#1}}
\newcommand{\mathfrak}[1]{\frak{#1}}

\newcounter{mycomment}
\newcommand\Hilb{\mathcal{H}}

\newcommand\RR{\mathbb{R}}
\newcommand\slim{\operatorname{s-}\lim}
\newcommand\ran{\operatorname{ran}}
\newcommand\Sch{\mathfrak{S}}
\newcommand\sgn{\operatorname{sgn}}
\newcommand\Imag{\operatorname{Im}}
\newcommand\Id{\operatorname{Id}}
\def\ang#1{\langle #1 \rangle}
\newcommand\sql{\sqrt{\lambda}}

\newtheorem{thm}{{\sc Theorem}}[section]

\newtheorem{lem}[thm]{{\sc Lemma}}

\newtheorem{prop}[thm]{{\sc Proposition}}
\newtheorem{defn}[thm]{{\sc Definition}}

\numberwithin{equation}{section}

\begin{document}

\title[Lectures on scattering theory] {Lectures on scattering theory}

\author{Dmitri Yafaev}
\address{Department of Mathematics, University of Rennes -- I, Campus de Beaulieu 35042 Rennes,
FRANCE}
\email{yafaev@@univ-rennes1.fr}

\maketitle
\centerline{\it Lecture notes prepared by Andrew Hassell, based on lectures }
\centerline{\it given by the author at the Australian National University }
\centerline{\it  in October and November, $2001$.}

\bigskip\bigskip

The first two lectures are devoted to describing the basic concepts of scattering theory
 in a very compressed way. A
detailed presentation of the abstract part can be found in \cite{I} and numerous applications in
\cite{RS} and \cite{Y2}. The last two lectures are based on the recent research of the author.

\bigskip

 \section{Introduction to Scattering theory. Trace class method}
 
{\bf 1}. Let us, first, indicate the place of scattering theory amongst other mathematical
theories.  This is very simple: it is a subset of perturbation theory. 
The ideology of perturbation theory is as follows.
Let $H_0$ and $H= H_0 + V$ be self-adjoint operators on a Hilbert space $\Hilb$, and  let $V$ be,
 in some sense, small compared to $H_0$. Then it is expected that the spectral properties of $H $
are close to those of $H_0$. Typically, $H_0$ is simpler than $H$, and in many cases we know 
its spectral family $E_0(\cdot)$   explicitly. 
The task of perturbation theory is to deduce information about the spectral properties of $H = H_0
+ V$ from those of $H_0$.
We shall always consider the case of self-adjoint operators on a Hilbert space $\Hilb$. 

   The spectrum of a self-adjoint operator has two components: discrete (i.e., eigenvalues) and
continuous. Hence,   perturbation theory has two parts: perturbation theory for the discrete
spectrum, and for the continuous spectrum. Eigenvalues of $H_0$ can generically be shifted under
arbitrary small perturbations, but the formulas  for these shifts are basically the same as in the
finite dimensional case, $\dim \Hilb < \infty$, which is   linear algebra.  On the contrary, the
continuous spectrum is much more stable.
 
\medskip

{\it Example.} Let 
\begin{equation}
H_0 = -\Delta,\quad H = -\Delta + v(x),\quad v(x)=\overline{v(x)}, \quad v\in L^\infty (\RR^d),
\label{Schr}\end{equation}
be self-adjoint operators on the Hilbert space
$\Hilb = L^2(\RR^d)$.  Denote by $\Phi$   the Fourier transform. Then 
\begin{equation} 
H_0 = \Phi^* |\xi|^2 \Phi,
\label{EH}\end{equation} 
 so that the spectrum of
$H_0$ is the same as that  of the multiplication operator $|\xi|^2$. Moreover, we know all
functions of the operator $H_0$ and, in particular, its spectral family explicitly.  In the case
considered $V$ is multiplication by  the function $v$.
 Thus,   each of the operators $H_0$ and $V$ is very simple. However, it is not
quite so simple to understand their sum, the Schr\"odinger operator $H = H_0+V$. Nevertheless, it
is easy to deduce from the Weyl theorem that, if $v(x)\rightarrow 0$ as $|x|\rightarrow\infty$,
then the essential spectrum of
$H$ is the same as that of $H_0$ and hence it coincides with $[0,\infty)$.

\medskip 

  Scattering theory requires classification of the spectrum in terms of the theory of measure. Each
measure may be decomposed into three parts: an absolutely continuous part, a singular continuous
part, and a pure point part. The same classification is valid for the spectral measure $E(\cdot)$
of a self-adjoint operator $H$. Thus, there is a decomposition of the Hilbert space $\Hilb =
\Hilb^{ac}\oplus\Hilb^{sc} \oplus\Hilb^{pp}$ into the orthogonal sum of invariant subspaces of the
operator
$H$; the operator restricted to $\Hilb^{ac}$,
$\Hilb^{sc}$ or $\Hilb^{pp}$ shall be denoted $H^{ac}$, $H^{sc}$ or $H^{pp}$, respectively. 
 The pure point   part corresponds to eigenvalues. The singular continuous part is typically
absent.  Actually, a part of scattering theory is devoted to proving this for various operators of
interest, but in these lectures we study only the  absolutely continuous part $H^{ac}$ of $H$.
The same objects for the operator $H_0$ will be labelled by the index `$0$'. 
We denote $P_0$   the orthogonal projection onto the absolutely continuous subspace $\Hilb_0^{ac}$
 of $H_0$.

The starting point of scattering theory is that the absolutely continuous part of self-adjoint
operator  is stable under fairly general perturbations. However assumptions on perturbations are
much more restrictive than those required for stability of the essential spectrum. So scattering
theory can be defined as perturbation theory for the absolutely continuous spectrum. Of course, it
is too much to expect that $H^{ac} = H_0^{ac}$. However, we can hope for a unitary equivalence:
$$
H^{ac} = U H_0^{ac} U^*, \quad U : \Hilb_0^{ac} \to \Hilb^{ac} \; \text{onto}. 
$$
The first task of scattering theory is to show this unitary equivalence. 
We now ask: how does one find such a unitary equivalence $U$? This is related (although the
relationship is not at all obvious) to the second task of scattering theory, namely, the large time
behaviour of solutions
$$ 
u(t) = e^{-iH t} f.
$$
 of the time-dependent equation
$$
i \frac{ \partial u}{\partial t} = H u, \quad u(0) = f \in \Hilb.
$$
If $f$ is an eigenvector, $H f = \lambda f$, then   $u(t) = e^{-i \lambda t} f$, so the   time
behaviour is evident. By contrast, if $f \in \Hilb^{ac}$, one cannot, in general, calculate $u(t)$
explicitly, but scattering  theory allows us to find its asymptotics as $t\rightarrow\pm\infty$.
 In the perturbation  theory setting,
it is natural to understand the asymptotics of $u$ in terms of solutions of the unperturbed equation,
$i u_t = H_0 u$. It turns out that, under rather general assumptions, for all $f \in \Hilb^{ac}$,
there are
$f_0^\pm\in\Hilb_0^{ac}$ such that
$$
u(t) \sim u_0^\pm(t),\  t \to \pm \infty, \text{ where } u_0^\pm(t) = e^{-i H_0 t} f_0^\pm,
$$
or, to put it differently,
\begin{equation}
\lim_{t \to \pm \infty}\big\| e^{-i H t} f - e^{-i H_0 t} f_0^\pm \big\| = 0.
\label{WOP}\end{equation}
Hence   $f_0^\pm$ and $f$ are related by the equality   
$$
f = \lim_{t \to \pm \infty} e^{i H t} e^{-i H_0 t} f_0^\pm,
$$
which justifies the following fundamental definition  given by C. M\o ller
\cite{Mo} and   made precise by K. Friedrichs \cite{Fr1}.

 \begin{defn}\label{WO}
 The limit
$$
W_\pm= W_\pm(H, H_0) = \slim_{t \to \pm \infty} e^{i H t} e^{-i H_0 t} P_0,  
$$
if it exists, is called the wave operator.
 \end{defn}

It follows that $f = W_\pm f_0^\pm$. The wave operator has the properties

(i) $W_\pm$ is isometric on $\Hilb_0^{ac}$. 

(ii) $W_\pm H_0 = H W_\pm$ (the {\it intertwining property}).

In particular, $\Hilb_0^{ac}$ is unitarily equivalent, via $W_\pm$, to the restriction 
$H |_{\ran W_\pm}$ of $H$ on the range $\ran W_\pm$ of the wave operator $W_\pm$ and hence
$\ran W_\pm\subset\Hilb^{ac}$. 

 \begin{defn}
If $\ran W_\pm=\Hilb^{ac}$, then $W_\pm$ 
is said to be complete. 
 \end{defn}

It is a simple result that $W_\pm(H, H_0)$ is complete if and only if the `inverse' wave
operator $W_\pm(H_0, H)$ exists.  Thus, if $W_\pm$ exists and is complete (at least for one of the
signs),
$H_0^{ac}$ and
$H^{ac}$ are unitarily equivalent. It should be emphasized that scattering theory is interested 
only in the canonical unitary equivalence provided by the wave operators. 

Another important object of scattering theory 
\begin{equation}
 {\bold S} = W_+^* W_-  
\label{SOP}\end{equation} 
is called the scattering operator.  
It commutes with $H_0$,
${\bold S} H_0 = H_0 {\bold S}$, which follows directly from property (ii) of the wave operators.
Moreover, it is unitary on $\Hilb^{ac}_0$ if both $W_\pm$ are complete. In the spectral
representation of the operator $H_0$, the operator ${\bold S}$ acts as multiplication by the
operator-valued function
$S(\lambda)$ known as the scattering matrix (see the next lecture, for more details).

An important generalization of Definition~\ref{WO} is due to   Kato \cite{Ka3}.

 \begin{defn} 
Let $J$ be a bounded operator. Then the modified wave operator 
$W_\pm(H, H_0, J)$ is defined by
\begin{equation}
 W_\pm(H, H_0, J) = \slim_{t \to \pm \infty} e^{i H t} J e^{-i H_0 t}P_0,
\label{WOJ}\end{equation}
 when this limit exists.
 \end{defn}
 
Modified wave operators still enjoy the intertwining property
$$ 
W_\pm(H, H_0, J) H_0 = H W_\pm(H, H_0, J),
$$ 
 but of course their isometricity on
$\Hilb^{ac}_0$ can be lost.

\

{\bf 2}. 
We have seen that the wave operators give non-trivial spectral information about $H$. 
Thus, it is an important problem to find conditions guaranteeing the existence of wave operators.
There are two quite different methods: the trace class method, and the smooth method (see the next
lecture). The trace class method is the principal method of abstract scattering theory.
 For applications to differential operators, both methods are important. 

The fundamental theorem for the trace class method is the Kato-Rosenblum theorem 
\cite{Ka1,R,Ka2}.  Recall
  that a compact operator $T$  on $\Hilb$ is in the class $\Sch_p$, $p>0$, if 
$$
|| T||_p^p = \sum \big(\lambda_j (T^* T)\big)^{p/2} < \infty.  
$$
In particular,   $\Sch_1$ is called the trace class and $\Sch_2$ is called the Hilbert-Schmidt
 class.

\begin{thm}\label{KR}
 If the difference $V = H - H_0$ belongs to the trace class, then
the wave operators $W_\pm(H, H_0)$ exist. 
\end{thm}

This is a beautiful theorem. It has a number of advantages, including:

(i) Since the conditions are symmetric with respect to the operators $H_0$ and $H$, the   wave
operators $W_\pm(H_0, H)$ also exist and hence  $W_\pm(H, H_0)$ are complete.

(ii) The formulation is simple, but all proofs of it are rather complicated. 

(iii) It relates very different sorts of mathematical objects: operator ideals, and scattering
theory.

(iv) It is effective, since it is usually easy to determine whether $V$ is trace class. 

(v) It is sharp, in the sense that if $H_0$ and $p > 1$ are given, 
there is a $V \in \Sch_p$ such that the spectrum of $H_0 + V$ is purely point.

However, Theorem~\ref{KR} has a disadvantage: it is useless in applications to differential
operators. Indeed, for example, for the pair \eqref{Schr}   $V$ is a multiplication operator which
cannot be even compact (unless identically zero), and therefore Theorem~\ref{KR} does not work.
Nevertheless, it is still useful if one is only interested in the unitary equivalence of
$H^{ac}$ and $H_0=H_0^{ac}$. In fact, it is sufficient to show that the operators
$  (H + c)^{-1}$ and
$ (H_0 + c)^{-1}$ are unitarily equivalent for some $c>0$ or, according to Theorem~\ref{KR}, that
their difference is trace class.

Actually,   this condition   is sufficient for the
existence (and completeness) of the wave operators. More generally, the following result is true.

\begin{thm}\label{Bir1} 
  Suppose that
 \[
 (H-z)^{-n}-(H_0-z)^{-n}\in{\goth S}_1
\] 
for some $n=1,2,\ldots$ and all $z$ with $\Imag z\neq 0$. Then the wave operators $W_\pm (H,H_0)$
exist and are complete.
\end{thm}

 Theorem~\ref{Bir1} was proved in \cite{BK} for $n=1$ and in \cite{Ka4} for arbitrary $n$.
  For semibounded operators, Theorem~\ref{Bir1} follows from
the Invariance Principle, due to Birman \cite{Birman}.

\begin{thm}\label{inv} 
 Suppose that
$\varphi(H)- \varphi(H_0)\in{\goth S}_1$ for a real function $\varphi$ such that its derivative
$\varphi^\prime$ is absolutely continuous and $\varphi^\prime(\lambda)> 0$. Then the wave operators
$W_\pm (H,H_0)$ exist and
\[
 W_\pm (H,H_0)=W_\pm (\varphi(H), \varphi(H_0)).
\]
The operators $W_\pm(H, H_0)$ should be replaced here by $W_\mp(H, H_0)$ if $\varphi'$ is negative.
\end{thm}

\bigskip

 {\bf 3}.
A typical result of trace class theory in applications to differential operators is the
following 

\begin{thm}\label{ScrTr}
  Let  
$H_0 = -\Delta+v_0(x)$  and $H = -\Delta +v_0(x)+ v(x)$ on $L^2(\RR^d)$. Assume that 
$v_0=\bar{v_0}$ is bounded and $v =\bar{v}$ satisfies
\begin{equation}
 | v(x) | \leq C \ang{x}^{-\rho},  \quad \ang{x} = (1 + |x|^2)^{1/2},
\label{power-decay}\end{equation}
 for some $\rho > d$.  Then the wave operators $W_\pm (H,H_0)$ exist and are complete.
\end{thm}

For simplicity we shall give the {\it proof}  only for $d = 1, 2$ or $3$.
We proceed from Theorem~\ref{Bir1} for $n=1$. Choose $c > 0$ large
enough so that $H_0 + c$ and
$H + c$ are invertible.   Then
 \[ 
(H + c)^{-1} - (H_0 + c)^{-1} = - (H + c)^{-1} V (H_0 + c)^{-1}.
\]
Denote temporarily $H_{00}=-\Delta$. Since the operators 
\[
(H + c)^{-1}(H_{00}+c)= \Id - (H + c)^{-1} (V_0+V)
\]
and, similarly, $(H_{00} + c) (H_0 +c)^{-1}$ are bounded, it suffices to check that
 \begin{equation}
 (H_{00} + c)^{-1} V (H_{00} + c)^{-1}= \Big( (H_{00} + c)^{-1} |V|^{1/2} \Big) \sgn v \Big(
|V|^{1/2}  (H_{00} + c)^{-1}  \Big)\in{\goth S}_1.
\label{HH}\end{equation} 
Note that $\Phi(H_{00} + c)^{-1} |V|^{1/2}$ is an integral operator with  kernel
$$
 (2\pi)^{-d/2} e^{-ix \cdot \xi} (|\xi|^2 + c)^{-1} |v(x)|^{1/2},
$$ 
 which is evidently in $L^2(\RR^{2d})$. Thus,
the operators  $(H_{00} + c)^{-1} |V|^{1/2}$  and its adjoint $|V|^{1/2} (H_{00} +
c)^{-1} $   are Hilbert-Schmidt and hence \eqref{HH} holds.   

Using Theorem~\ref{Bir1} for  
$n>d/2-1$, it is easy (see \cite{BY}, for details) to extend this result    to an arbitrary
$d$. On the contrary, the condition    $\rho >d$ in
 \eqref{power-decay}  cannot be improved in the trace-class framework.

\bigskip

 {\bf 4}.
 Theorem~\ref{KR} admits the following generalization (see \cite{Pe})  to the wave operators
\eqref{WOJ}.

\begin{thm}\label{P}
 Suppose that   $V=HJ-JH_0\in{\goth S}_1$. Then the
wave operators $W_\pm (H,H_0;J)$ exist.
\end{thm}

This result due to Pearson allows to simplify considerably the original proof of
Theorem~\ref{KR}  and of its different generalizations.
Although still rather sophisticated, the proof of Theorem~\ref{P} relies only on the following 
elementary lemma of   Rosenblum.

\begin{lem}\label{R}
 For  a self-adjoint operator $H$, consider the set   ${\goth R} \subset {\cal
H}^{ac}$ of elements
$f$ such that
\[ r_H^2(f):=\text{ess-sup}\: d(E (\lambda)f,f)/d\lambda <\infty.
\]
 If $G $   is a Hilbert-Schmidt
operator, then for any $f\in{\goth R} $ 
\[
\int_{-\infty}^\infty ||G\exp(-iHt)f||^2dt \leq 2\pi r^2_H (f) ||G||^2_2.
\]
 Moreover, the set ${\goth R} $ is dense in   ${\cal H}^{ac}$.
\end{lem} 

   The wave operators were defined in terms of exponentials, $e^{-i t H}$.
 There is also a `stationary' approach, in which the exponentials are replaced   by the resolvents
 $R(z)=(H -z)^{-1}$, $z=\lambda\pm i\varepsilon$,
and the limit   $t\rightarrow\pm\infty$ is replaced by the limit 
$\varepsilon\rightarrow 0$.  In the trace class framework a consistent stationary approach was
developed  in the paper
\cite{BE}. From analytical point of view it relies on the following result 
on boundary values of  
resolvents which is interesting in its own sake.

\begin{prop}\label{BE} 
 Let $H$ be a  self-adjoint operator and let $G_1$, $G_2$ be 
arbitrary Hilbert-Schmidt operators. Then the operator-function
$G_1 R(\lambda+i\varepsilon)  G_2$ has   limits as $\varepsilon\rightarrow 0$ $($and $G_1
R(z)  G_2$ has angular boundary values as $z\rightarrow \lambda\pm i0)$ in the Hilbert-Schmidt
class for almost all
$\lambda\in{\Bbb R}$. Moreover, the operator-function
$G_1 E(\lambda )  G_2$ is differentiable in the trace norm for almost all
$\lambda\in{\Bbb R}$.
\end{prop}

In particular, Proposition~\ref{BE} allows one to obtain a stationary proof of the Pearson theorem
(see   the book \cite{I}).

\section{Smooth method. Short range scattering theory}

 {\bf 1}. 
The   smooth  method relies on a certain regularity of the perturbation in the spectral 
representation of the operator $H_0$. 
 There are different ways to understand  regularity. For example, in the Friedrichs-Faddeev model
\cite{F} $H_0$ acts as multiplication by independent variable in the space ${\cal H}=L_2(\Lambda;
{\goth N})$ where $\Lambda$ is an interval  
and ${\goth N}$ is an auxiliary Hilbert space. The perturbation
$V$ is an integral operator with sufficiently smooth kernel.

 Another possibility is to use the concept of $H$-smoothness introduced by T. Kato in \cite{Ka5}.

\begin{defn} Let $G$ be an $H$-bounded operator;  that is, suppose that $G (H + i)^{-1}$ is
bounded. Then we say that $G$ is $H$-smooth if there is a $C < \infty$ such that
\begin{equation}
\int_{-\infty}^\infty \| G e^{-i H t} f \|^2 dt \leq C^2 \| f \|^2 
\label{smooth-cond}\end{equation}
 for all $f \in \Hilb$ or, equivalently,  
\begin{equation}
\sup_{\Imag z \neq 0} \| G \big( R(z) - R(\overline{z})\big) G^* \| \leq 2\pi C^2.    
\label{smooth-condA}\end{equation} 
\end{defn}

  In applications the assumption of $H$-smoothness of an operator $G$ imposes too stringent
conditions on the operator $H$. In particular, the operator $H$ is necessarily absolutely
continuous if      kernel of $G$  is trivial.  This excludes eigenvalues   and other
singular points in the spectrum of
$H$, for  example, the bottom of the continuous spectrum  for the Schr\"odinger operator with
decaying potential or edges of bands if the spectrum has the band structure. However it is often
suffices to verify   $H$-smoothness of the operators $G E(X_n)$ where the union of 
intervals $X_n$ exhausts ${\Bbb R}$ up to a set of the Lebesgue measure zero. In this
case we say that $G$ is locally $H$-smooth.

We have the following

\begin{thm}\label{smooth} 
Let $H - H_0 = G^* G_0$, where
$G_0 $ is locally $H_0$-smooth and $G $ is locally
$H$-smooth. Then the wave operators $W_\pm(H, H_0)$ exist and are complete.
\end{thm}

This is a very useful theorem, yet the proof is totally elementary. To make it even more simple,
we forget about the word `locally' and assume that
 $G_0  $ is $H_0$-smooth and $G  $ is $H$-smooth. We write
$$
\lim_{t \to \pm\infty} \big( e^{iH t} e^{-t H_0 t} f_0, f \big) = \big(  f_0, f \big)+
\lim_{t \to \pm\infty} i \int_0^t \big( G_0 e^{-s H_0 } f_0, G e^{-i H s} f \big) \, ds.
$$ Since the left and right hand sides of the last inner product are $L^2({\Bbb R})$,
 the limit on the right hand side exists. This shows existence of the {\it weak} limit. To show the
strong limit, we estimate using \eqref{smooth-cond}
$$\begin{gathered}
\big|\big( e^{iH t} e^{-t H_0 t} f_0, f \big) - \big( e^{iH t'} e^{-t H_0 t'} f_0, f \big)\big|
 =  
\big|\int_{t'}^t \big( G_0 e^{-s H_0 } f_0, G e^{-i H s} f \big) \, ds \big|
\\
\leq C \| f \|  \Bigl(\int_{t'}^t \| G_0 e^{-s H_0 } f_0 \|^2 \, ds\Bigr)^{1/2}.
\end{gathered}$$
 Taking the sup over $f$ with $\| f \| = 1$, we obtain
$$
\| e^{iH t} e^{-t H_0 t} f_0 - e^{iH t'} e^{-t H_0 t'} f_0 \| \leq C  \Bigl(\int_{t'}^t \| G_0
e^{-s H_0 } f_0 \|^2 \, ds \Bigr)^{1/2}
$$ which goes to zero as $t', t$ tend to infinity.

Of course Theorem~\ref{smooth} is not effective since the verification of $H_0\,$- and especially
of $H$-smoothness may be a difficult problem. In the following assertion the hypothesis only
concerns the free resolvent, $R_0(z)=(H_0-z)^{-1}$.

\begin{thm}\label{GR_0} Suppose $V$ can be written in the form $V = G^* \Omega G$,
 where $\Omega$ is a bounded operator and $G R_0(z) G^*$ is compact for $\Imag z \neq 0$, and is
norm-continuous up to the real axis except possibly at a finite number of points $\lambda_k$,
$k=1,\ldots, N$.    Then
$W_\pm(H, H_0)$ exist and are complete.
\end{thm}

 We  give only a brief sketch of its proof. Set
\[
X_n=(-n,n)\setminus \bigcup_{k=1}^N (a_k-n^{-1},a_k+n^{-1}).
\]
Since $\| G R_0(\lambda\pm i\varepsilon) G^* \|$ is uniformly bounded for $\lambda\in
X_n$,  the operator   $GE_0(X_n)$ is $H_0$-smooth (cf.   the definition (\ref{smooth-condA})).

To show a similar result for the operator-function   $G R (z) G^*$,  we use the resolvent identity
\begin{equation}
 R(z) = R_0(z) - R_0(z) V R(z),
\label{Res}\end{equation}
 whence
\begin{equation}
 G R(z) G^* = \Big( \Id + G R_0(z) G^* \Omega \Big)^{-1} G R_0(z) G^*,\quad \Imag z\neq 0.
\label{GRG}\end{equation}
It easily follows from self-adjointness of the operator $H$ that the homogeneous equation  
\begin{equation}
f + G R_0(z) G^* \Omega f=0
\label{hom}\end{equation}
has only the trivial solution $f=0$ for
$\Imag z\neq 0$. Since the operator $G R_0(z) G^*$ is compact, this implies
that  the inverse operator in \eqref{GRG} exists. By virtue of equation \eqref{GRG}, the
operator-function   $G R (z) G^*$ is continuous up to the real axis except the points $\lambda_k$
and the set
${\cal N}$ of
$\lambda$ such that the equation \eqref{hom} has a non-trivial solution for $z=\lambda+i0$ or
$z=\lambda-i0$. The set
${\cal N}$ is obviously closed. Moreover, it has the Lebesgue measure zero  by the analytical
Fredholm alternative. This implies that the pair $G$, $H$ also satisfies the conditions of
Theorem~\ref{smooth}. It remains to use this theorem.

\

{\bf 2}.
Let us return  to the Schr\"odinger operator $H = -\Delta + v$.
Potentials $v$ satisfying \eqref{power-decay} with $\rho > 1$ are
 said to be {\it short range}. Below in this lecture we make this assumption. Let us apply
Theorem~\ref{GR_0} to the pair $H_0=-\Delta$, $H$. Put $r=\rho/2$ and  $G = \ang{x}^{-r}$. One
can verify (cf. the proof of Theorem~\ref{ScrTr}) that  
$G R_0(z) G^*$ is compact for $\Imag z \neq 0$ (and arbitrary $r >0$). Next we consider the
spectral family
$E_0(\lambda)$ which according to \eqref{EH} satisfies
\begin{equation} 
d(E_0(\lambda)f, f)/d\lambda = ||\Gamma_0(\lambda) f||^2,
\label{sob}\end{equation}
where
\begin{equation}
(\Gamma_0(\lambda) f) (\omega) = 2^{-1/2} \lambda^{(d-2)/4} \hat f(\sqrt{\lambda} \omega) 
\label{Gamma}\end{equation}
and $\hat f=\Phi f$ is the Fourier transform of a function $f$ from, say,
the Schwartz class. Thus, up to a numerical factor, $\Gamma_0(\lambda) f$ is the restriction of
$\hat f$ to the sphere of radius  $\sqrt{\lambda}$.
 Remark further that if $f \in \ang{x}^{-r} L^2(\RR^d)$, then $\hat f$
 belongs to the Sobolev space $\text{H}^{r}(\RR^d)$. Since $r  > 1/2 $, it follows from  the
Sobolev trace theorem that the operator
\[
\Gamma_0(\lambda) \ang{x}^{-r} : L^2(\RR^d)\rightarrow L^2({\Bbb S}^{d-1})
\]
is bounded and  
  depends (in the operator norm) H\"older continuously on $\lambda>0$. Therefore according to
\eqref{sob}
 the operator-function $GE_0(\lambda)G$ is differentiable and its derivative is also
   H\"older continuous. Now we use the
representation 
$$ 
(R_0(z) G f,G f) = \int_0^\infty (\lambda - z)^{-1} d(E_0(\lambda) Gf,G f) ,
$$
which, by the Privalov theorem, implies that the operator-function
 $ G R_0(z) G $
is continuous in the closed complex plane cut along $[0,\infty)$ except, possibly, the point
$z=0$.  Applying Theorem~\ref{GR_0}, we now obtain

\begin{thm}\label{SR} 
Let $v$ be short range.    Then the wave opeartors
$W_\pm(H, H_0)$ for the pair $\eqref{Schr}$ exist and are complete.
\end{thm}

 Theorem~\ref{SR} implies that for every $f \in \Hilb^{ac}$, there is $f_0^\pm$ such that
relation \eqref{WOP} holds. Using the well-known expression for kernel of the integral operator
$e^{-it H_0}  $ (in the $x$-representation), we find
that
\begin{equation}
 (e^{-it H} f)(x) \sim   e^{i|x|^2/4t} (2  i t)^{-d/2} \hat f_0^\pm(x/2t).
\label{SRX}\end{equation} 
Here `$\sim$' means that the difference of the left  and right hand sides tends to zero in
$L^2({\Bbb R}^d)$ as
$t\rightarrow\pm\infty$. Thus, the solution `lives' in the region $|x| \sim |t|$ of $(x, t)$
space. 

As a by-product of our considerations we obtain that the operator-function
$GR(z)G$ is continuous up to real axis, except     a  closed set of
measure zero. More detailed analysis shows that this set consists of eigenvalues of the operator
$H$ (and possibly the point zero), so that the singular continuous spectrum of $H$ is empty.
Finally, we note that, by the Kato theorem, the operator $H$ does not have positive eigenvalues.
This gives the following assertion known as the limiting absorption principle.

\begin{thm}\label{LAP}
  Let $v$ be short range and $r>1/2$.    Then the operator-function
 $ \ang{x}^{-r}R (z) \ang{x}^{-r} $ is norm-continuous in the closed complex plane cut along
$[0,\infty)$ except negative eigenvalues of the operator $H$ and, possibly, the point $z=0$. 
\end{thm}
   
 \

{\bf 3}.
Let us compare Theorems~\ref{ScrTr} and \ref{SR}. If $v_0=0$, then Theorem~\ref{SR}  is stronger
because, in assumption \eqref{power-decay}, it requires 
   that $\rho>1$   whereas Theorem~\ref{ScrTr} requires 
   that $\rho>d$. Theorem~\ref{SR} can be extended to some
other cases, for example to periodic and long-range        $v_0$. In the first case the
spectral family   $E_0(\cdot)$ can
be constructed rather explicitly. In the second case the limiting absorption principle can be also
verified (see the next lecture). However, contrary to Theorem~\ref{ScrTr}, the method of proof of
Theorem~\ref{SR} gives nothing for arbitrary $v_0
\in L^\infty$. Therefore the following question naturally arises.

\medskip

{\it Problem}.  Let  
$H_0 = -\Delta+v_0(x)$  and $H = -\Delta +v_0(x)+ v(x)$ where $v$ satisfies estimate
\eqref{power-decay} for $\rho>1$. Do the wave operators
$W_\pm(H, H_0)$ exist for an arbitrary $v_0 \in L^\infty$?

\medskip 

This problem is of course related to unification of trace class and smooth approaches.
  The following theorem does this to some extent:

\begin{thm}\label{G-bvs}
 Assume that $V$ is of the form $V = G^* \Omega G$, where $\Omega $ is a
bounded operator, $G R_0(z) G^* \in \Sch_p$ for some $p < \infty$, and $G  R_0(z) G^*$ has angular
boundary values in $\Sch_p$ for almost every $\lambda \in \RR$. Then $W_\pm(H, H_0)$ exist and are
complete. 
\end{thm}

Note that, given  Proposition~\ref{BE}, Theorem~\ref{G-bvs} provides an independent proof of
Theorem~\ref{KR}. On the other hand, it resembles Theorem~\ref{GR_0} of the smooth approach.
However it gives nothing for the solution of the problem formulated above. Most probably, the
answer to the formulated question  is negative, which can be considered as a strong evidence that
a real  unification of trace class and smooth approaches does not exist.

 \

{\bf 4}.
 The stationary method is intimately related to eigenfunction expansions of the operators $H_0$
and $H$. Let us discuss this relation on the example of
the Schr\"odinger operator $H$ with a short range potential. For the operator $H_0=-\Delta$, the
construction of eigenfunctions is obvious.
Actually, if   
$\psi_0(x,   \omega,\lambda)=e^{i \sqrt{\lambda} \omega \cdot x}$,
 then
$ -\Delta \psi_0 = \lambda \psi_0$. This collection of eigenfunctions is `complete',
 so eigenfunctions are parametrized by $\omega \in {\Bbb S}^{d-1}$ (for fixed $\lambda > 0$). 
By the intertwining property, the wave operators $W_\pm(H, H_0)$  diagonalize $H$ and hence
$$
  H W_\pm \Phi^* = |\xi|^2 W_\pm \Phi^*.
$$
 Thus, at least formally, eigenfunctions of $H$, that is, solutions of the equation
\begin{equation}
  -\Delta \psi + v \psi = \lambda \psi,
\label{efn-eqn}\end{equation}
can be constructed by one of the two equalities $\psi_+(  \omega,\lambda)=W_+ \psi_0( 
\omega,\lambda)$ or $\psi_-(  \omega,\lambda)=W_- \psi_0( \omega,\lambda)$.

  It turns out that this definition can be given a precise sense, and one can construct solutions
of the Schr\"odinger equation  with   asymptotics $\psi_0(x,   \omega,\lambda)$ at infinity. 

\begin{thm}
 Assume \eqref{power-decay} is valid for some $\rho > d$. For every $\lambda > 0$ and $\omega \in
{\Bbb S}^{d-1}$ there is a solution $\psi$
 of \eqref{efn-eqn} such that
\begin{equation}
\psi(x,   \omega,\lambda) = e^{i \sqrt{\lambda} \omega \cdot x} + a(\hat x,  \omega,\lambda)
|x|^{-(d-1)/2} e^{i \sqrt{\lambda} |x|} + o(|x|^{-(d-1)/2}), \quad \hat x=x/|x|,
\label{efn-exp}\end{equation}
where $a$ is a continuous function on ${\Bbb S}^{d-1}\times {\Bbb S}^{d-1}$.
\end{thm}

We interpret the $e^{i \sqrt{\lambda} \omega \cdot x}$ term as an incoming plane wave, 
the $|x|^{-(d-1)/2} e^{i \sqrt{\lambda} |x|}$ term as an outgoing spherical wave, and the
coefficient $a(\hat x,\omega,\lambda )$ is called the {\it scattering amplitude} for the incident
direction $\omega$ and the direction of  observation $\hat x$. Note that  a solution  of the
Schr\"odinger equation is determined uniquely by the condition that, asymptotically, it is a sum
of a  plane and the outgoing spherical waves.

Under assumption \eqref{power-decay} where $\rho > (d+1)/2$, eigenfunctions of the operator $H$
 may be constructed by means of the following formula:
\begin{equation}
\psi (  \omega,\lambda) = \psi_0(  \omega,\lambda) - R(\lambda + i0)  V
\psi_0(\omega,\lambda) ,
\label{psi}\end{equation}
or, equivalently, as solutions of the Lippman-Schwinger equation
$$
\psi (  \omega,\lambda) = \psi_0(  \omega,\lambda) + R_0(\lambda + i0)  V \psi( 
\omega,\lambda).
$$
We note that, strictly speaking, the second term in the right hand side of \eqref{psi} is 
defined by the equality  
\begin{equation}
 R(\lambda + i0)  V\psi_0(\omega,\lambda)=
  R(\lambda + i0) \ang{x}^{-1/2 - \epsilon}    (V \ang{x}^{(d+1)/2 + 2\epsilon}) (
\ang{x}^{-d/2 - \epsilon} \psi_0(  \omega,\lambda) ),
\label{corr}\end{equation}
where $\varepsilon>0$ is sufficiently small. Here  $\ang{x}^{-d/2 - \epsilon} \psi_0(\lambda,
\omega) \in L^2({\Bbb R}^d)$ and the operator $V \ang{x}^{(d+1)/2 + 2\epsilon}$ is bounded.
Therefore, by the limiting absorption principle (Theorem~\ref{LAP}), the function \eqref{corr}
belongs to the space $\ang{x}^{1/2 + \epsilon} L^2({\Bbb R}^d)$ for any $\epsilon>0$. The
asymptotics
\eqref{efn-exp} is still true for arbitrary $\rho > (d+1)/2$ if  the remainder
$e(x)=o (|x|^{-(d-1)/2})$ is understood in the following averaged sense: 
$$
 e(x) = o_{av}(|x|^{-(d-1)/2}) \Longleftrightarrow \lim_{R \to \infty} \frac1{R} \int_{|x| \leq R}
|e(x)|^2 dx = 0.
$$
In this case the scattering amplitude $a(\hat x,  \omega,\lambda)$ belongs to the space
$L^2({\Bbb S}^{d-1})$ in the variable $\hat x$ uniformly in $\omega\in{\Bbb S}^{d-1}$.

 It is often convenient
 to write $\psi$ in terms of the parameter  
$\xi =\sqrt{\lambda}\omega \in\RR^d$, instead of $(\omega,\lambda )$. Thus, we set
$$
\psi_-(x, \xi) = \psi(x,   \omega,\lambda), \quad
\psi_+(x, \xi) = \overline{\psi_-(x, -\xi)}.
$$
The Schwartz kernels of the wave operators are intimately related to eigenfunctions.
 In fact, if we define `distorted Fourier transforms' $\Phi_\pm$  by
$$
(\Phi_\pm f)(\xi) = (2\pi)^{-d/2} \int_{{\Bbb R}^d} \overline{\psi_\pm(x, \xi)} f(x) dx,
$$
so that $\Phi_\pm H = |\xi|^2 \Phi_\pm$, then we have
\begin{equation}
W_\pm = \Phi^*_\pm \Phi_0,
\label{Phi}\end{equation}
where  $\Phi_0=\Phi$ is the classical Fourier transform.
Notice that the right hand side is defined purely in terms of time-independent quantities. 
This can be taken to be the {\it definition} of the wave operators
 in the stationary approach to scattering theory. Historically, it was the first approach to the
study of wave operators suggested by Povzner in \cite{Po1,Po2} and developed further in \cite{Ik}.
Under optimal assumption $\rho>1$ in \eqref{power-decay}  Theorem~\ref{SR} was obtained in
\cite{Ka7}.

It follows from  \eqref{Phi} that the scattering operator
\begin{equation}
{\bold S} = W_+^* W_-=\Phi^*_0 \Phi_+ \Phi^*_- \Phi_0.
\label{eq:PhiS}\end{equation}
Let $H_0$ be realized (via the Fourier transform) as multiplication by $\lambda$ in the space
$ L^2(\RR^+; L^2({\Bbb S}^{d-1}))$.
Since ${\bold S}$      commutes with $H_0$ and is unitary, it acts in this representation as
multiplication by the operator-function (scattering matrix)
$S(\lambda) : L^2({\Bbb S}^{d-1}) \to L^2({\Bbb S}^{d-1})$ which is also unitary for all
$\lambda>0$. It can be deduced from \eqref{eq:PhiS} that
\begin{equation}
 (S(\lambda) f)(\omega) = f(\omega) + \gamma c(\lambda)
\int_{{\Bbb S}^{d-1}} a(\omega, \omega',\lambda) f(\omega') d\omega',
\label{SAS}\end{equation}
 where $a$ is the scattering amplitude defined by \eqref{efn-exp} and 
\[
\gamma=e^{ \pi i(d-3)/4},\quad c(\lambda)=i   (2\pi)^{ -(d-1)/2}\lambda^{ (d-1)/4}.
\]

If assumption \eqref{power-decay} is satisfied for $\rho > 1$ only, then the scattering matrix
satisfies the relation
\begin{equation}
 S(\lambda) = \Id - 2\pi i \Gamma_0(\lambda) \Big( V - V R(\lambda + i0) V \Big)
\Gamma_0^*(\lambda),
\label{sc-stat}\end{equation}
which  generalizes \eqref{SAS}.
 Here  $\Gamma_0(\lambda)$ is the operator \eqref{Gamma}.  Since $\ang{x}^{-r} V \ang{x}^{-r}$ is
a bounded operator for $r=\rho/2>1/2$, we see that \eqref{sc-stat} is
correctly defined. It follows from formula \eqref{sc-stat} that the operator $S(\lambda) - \Id$ is
compact.
 Thus, the spectrum of $S(\lambda)$ (which lies on the unit circle by unitarity) is discrete, and
may accumulate only at the point $1$.   Moreover, if
 \begin{equation}
\big|  \partial^\alpha v(x) \big|  \leq C_{\alpha} \ang{x}^{-\rho  - |\alpha|}
\label{LR}\end{equation}
 for all multi-indices $\alpha$ (and $\rho>1$), then the
kernel  $k(\omega,
\omega')$ of the operator $S(\lambda) - \Id$ is smooth for
$\omega \neq
\omega'$  and   $|k(\omega,
\omega')| \leq C |\omega - \omega'|^{-d+\rho}$.

For $\rho\in ( 1,(d+1)/2]$,  construction of eigenfunctions, which behave asymptotically as plane
waves, becomes a difficult problem. In particular, formula \eqref{psi} makes no sense in this
case. One can do something, but  the construction of \cite{Skr}  is rather complicated and
requires the condition \eqref{LR}.
On the contrary,  for arbitrary $\rho > 1$ we can construct solutions which formally correspond to
averaging of
$\psi(x,\omega,\lambda)$ over
$\omega\in{\Bbb S}^{d-1}$. To illustrate this idea, let us first consider the free case $v=0$.
Then, for any $b
\in C^\infty({\Bbb S}^{d-1})$, the function
$$
 u(x) = \int_{{\Bbb S}^{d-1}} e^{i\sqrt{\lambda} \omega \cdot x} b(\omega) d\omega
$$ 
satisfies $-\Delta u = \lambda u$ and, by the stationary phase arguments, it has asymptotics
$$ 
u(x) = c(\lambda)^{-1} |x|^{-(d-1)/2}
\Big(\bar{\gamma} b(\hat x) e^{i \sqrt{\lambda} |x|} -  \gamma  b(-\hat x) e^{-i \sqrt{\lambda}
|x|}\Big) + o(|x|^{-(d-1)/2}).
$$
In the general case one can also construct   solutions of the
 Schr\"odinger equation \eqref{efn-eqn} with the asymptotics of incoming and outgoing spherical
waves.

\begin{thm}\label{aver}
 Assume \eqref{power-decay} is valid for some $\rho > 1$. 
Let $u$ be a solution of the Schr\"odinger equation \eqref{efn-eqn} satisfying
$$
\int_{|x| \leq R} |u(x)|^2 \, dx \leq C R, \text{ for all } R \geq 1.
$$
Then there are $b_\pm \in L^2({\Bbb S}^{d-1})$ such that
\begin{equation}
u(x) = |x|^{-(d-1)/2} \Big( \bar{\gamma} b_+(\hat x) e^{i \sqrt{\lambda} |x|}
-\gamma b_-(-\hat x) e^{-i
\sqrt{\lambda} |x|} \Big) + o_{av}(|x|^{-(d-1)/2}).
\label{sph-asympt}\end{equation}
Functions $b_\pm$ are related by the scattering matrix $:$ $b_+ =   S(\lambda) b_-$.
Conversely, for all
$b_+\in L^2({\Bbb S}^{d-1})$ $($or $b_-
\in L^2({\Bbb S}^{d-1}))$, there is a unique function $b_- \in L^2({\Bbb S}^{d-1})$
$($or $b_+\in L^2({\Bbb S}^{d-1}))$, and a
unique solution $u$ of
\eqref{efn-eqn} satisfying \eqref{sph-asympt}.
\end{thm}

\section{Long range  scattering theory}

There are different and, to a large extent, independent  methods  in long range  scattering (see 
\cite{Y2}). Here we shall give a brief presentation of the approach of the paper \cite{Y10} which
relies on the theory of smooth perturbations.

\medskip

{\bf 1}.
The condition \eqref{power-decay} with $\rho > 1$ is optimal even for the existence of wave
operators for the pair $H_0=-\Delta$, $H=-\Delta+v(x)$. For example,   the wave
operators          do not    exist if  $v(x)=v_0 \langle x\rangle^{-1}$, $v_0\neq 0$. Nevertheless
the asymptotic behaviour of  the function $ \exp (-iHt)f$ for large
$|t|$ remains sufficiently close to the free evolution $ \exp (-iH_0t)f_0$ if  the condition 
\eqref{LR} is satisfied for $  |\alpha|\leq\alpha_0$ with  $\alpha_0$   big
enough. Potentials obeying this condition for some
$\rho\in(0,1]$ are called long-range.

There are several possible descriptions of  $ \exp (-iHt)f  $ as  $t\rightarrow\pm\infty$. One of
them is a modification of the free evolution which, in its turn, can be done either in momentum
  or in coordinate  representations. Here we discuss the coordinate modification. Motivated by
\eqref{SRX}, we set
\begin{equation}
 (U_0 (t)f )(x) = \exp (i\Xi (x,t)) (2it)^{-d/2} \hat{f } (x/(2t)),  
\label{lr1}\end{equation}
where the choice of the phase function $\Xi$ depends on $v$. Then the wave operators are defined by
the equality
\begin{equation}
 W_\pm=   \slim_{t \to \pm \infty} e^{i H t} U_0 (t).  
\label{WOM}\end{equation}
To be more precise, these limits exist if the function $\Xi (x,t)$
 is  a  (perhaps, approximate) solution of the eikonal equation
\[ 
\partial \Xi/\partial t+ |\nabla \Xi|^2+ v=0.
\]
 For example, if $\rho>1/2$, we can neglect here the nonlinear term
$|\nabla\Xi|^2$   and  set
\begin{equation}
 \Xi(x, t) =(4t)^{-1} |x|^2  -t\int_0^1 v(sx) ds.
\label{WO6}\end{equation}
In the general case one obtains  an approximate  solution of the eikonal equation by
 the method of successive approximations.  
 With the phase $\Xi(x, t)$ constructed
in such a way,
   for an arbitrary $\hat{f}\in C_0^\infty ({\Bbb R}^d\setminus
\{0\})$,  the function $U_0 (t)f$ is an approximate solution of the time-dependent Schr\"odinger
equation in the sense that
\[  
\Bigl|\int_1 ^{\pm\infty} ||(i\partial/\partial t -H) U_0 (t)f|| dt \Bigr|<\infty.
\]
This condition implies that the vector-function $\partial e^{i H t} U_0 (t) f/\partial t\in
L^1({\Bbb R})$ and hence the limit \eqref{WOM} exists. The modified wave operators have all the
properties of usual wave operators. 
They are isometric,  $W_\pm H_0 = H W_\pm$ and 
$\ran  W_\pm \subset\Hilb^{ac}$.  As in the short-range case, the wave operator is said to be
complete if $\ran  W_\pm = \Hilb^{ac}$. Only the completeness of $W_\pm$ is a
non-trivial mathematical problem. As we shall see below, it has a positive solution which implies
that for every $f \in \Hilb^{ac}$, there is $f_0^\pm$ such that (cf. \eqref{SRX})
\begin{equation}
 (e^{-it H} f)(x) \sim   e^{i\Xi (x,t)} (2  i t)^{-d/2} \hat f_0^\pm(x/2t).
\label{LRX}\end{equation}
  Thus, in the short and long range cases, the large time asymptotics of  solutions of 
the time-dependent Schr\"odinger equation
differ only by the phase factor. In particular, in both cases they live  in the region $|x| \sim |t|$
of $(x, t)$ space. 

\bigskip

{\bf 2}.
For the proof  of completeness   of wave
operators, we need some analytical results which we discuss now. Note first of all that the
 limiting absorption principle, Theorem~\ref{LAP}, remains true if assumption \eqref{LR} is
satisfied for some $\rho>0$ and $  |\alpha|\leq 1$. However perturbative arguments of the previous
lecture do not work for long range potentials. The simplest proof
(see the original paper \cite{Mo1} or the book \cite{CFKS})  of  the
 limiting absorption principle relies  
in this case on the Mourre estimate:
\begin{equation}
 i E(I) [H, A] E(I) \geq c E(I), \quad A = -i \sum_{j=1}^d (\partial_j x_j
+ x_j \partial_j),\quad c>0,
\label{Mourre}\end{equation}
 where $I$ is a sufficiently small interval about a given $\lambda > 0$.
 This is not too hard to prove for long-range two body potentials (but is much harder for
multiparticle operators).  To be more precise, the Mourre estimate implies that  
the operator-function
 $ \ang{x}^{-r}R (z) \ang{x}^{-r} $ is norm-continuous up to positive half-axis, except a discrete
set of eigenvalues of the operator $H$. In particular, the operator $\ang{x}^{-r}$ is
locally $H$-smooth which is sufficient for scattering
theory. Moreover, independent arguments (see, e.g., \cite{RS}) show that the Schr\"odinger operator
$H$ does not  have positive eigenvalues. 
 
 Unfortunately, the operator $\ang{x}^{-1/2}$ is {\it not} (locally)
$H$-smooth even in the case $v=0$. However, there is a substitute:

\begin{thm}\label{G} 
Let assumption \eqref{LR} be satisfied for some $\rho>0$ and $  |\alpha|\leq 1$. Set
$$
(\nabla^\perp_j u)(x) = (\partial_j u)(x)  - |x|^{-2}((\nabla u)(x) \cdot   x)  x_j,
\quad j=1,\ldots, d.
$$ Then the operators $  \ang{x}^{-1/2} \nabla^\perp_j E(\Lambda)$ are $H$-smooth  for any compact $\Lambda
\subset (0,\infty)$. 
\end{thm}

The {\it proof } is based on the equality
$$ 
2|x|^{-1}\sum_{j=1}^d |\nabla^\perp_j u|^2=   ([H, \partial_r]u,u) + (v_1u,u),\quad
v_1(x)=  O(|x|^{- 1-\rho }),
$$ 
which is obtained by direct calculation. Since $\ang{x}^{-(1 + \rho)/2}E(\Lambda)$ is   $H$-smooth,
 we only have to consider the term $[H, \partial_r]$. For this we note that
$$ i \int_0^t \big( [H, \partial_r] e^{-i s H} f, e^{-i s H} f \big) \, ds
 = ( \partial_r e^{-i t H} f, e^{-i t H} f ) - (\partial_r f, f).
$$
 For $f \in \ran E(\Lambda)$,  
 the right hand side is bounded by $C(\Lambda) \| f \|^2$,  so this shows
  $H$-smoothness of the operators $ \ang{x}^{-1/2} \nabla^\perp_j E(\Lambda)$. 

Note that the commutator method used for the proof of Theorem~\ref{G} goes back to Putnam
\cite{Pu} and Kato \cite{Ka6}.

\bigskip

 {\bf 3}.
Our proof  of completeness    relies on consideration of wave opeartors
\eqref{WOJ} with a specially chosen operator $J$. Let us look for $J$ in the form of a
pseudodifferential operator
\begin{equation}
( Jf)(x) = (2\pi)^{-d/2} \int_{{\Bbb R}^d} e^{ix \cdot \xi} p(x, \xi) \hat f(\xi) d\xi.
\label{J}\end{equation}
We shall work in the class of symbols ${\cal S}^m_{\rho, \delta}$ consisting of functions 
$p(x, \xi)
\in C^\infty(\RR^d \times \RR^d)$ satisfying the estimates
$$
\big| \partial_x^\alpha \partial_\xi ^\beta p(x, \xi) \big| \leq C_{\alpha, \beta}
 \ang{x}^{m - \rho |\alpha| + \delta |\beta|} 
$$
for all multi-indices $\alpha$ and $\beta$. In addition, we shall assume that $p$ vanishes for
$|\xi|
\geq R$, for some $R$. The number
$m$ is called the {\it order} of the symbol, and of the corresponding pseudodifferential operator.  
Compared to the usual calculus, here the roles of $x$ and $\xi$ are interchanged: we require some
decay estimates as $|x| \to \infty$, rather than as $|\xi| \to \infty$ in the usual situation. We
shall assume that $ 0 \leq \delta < 1/2 < \rho \leq 1$. The symbol $p(x, \xi)$ of the operator
$J$ belongs to the class ${\cal S}^0_{\rho, \delta}$.
We recall that operators of order zero  are bounded on $L^2({\Bbb R}^d)$.

\medskip

Since  
$$
 de^{i H t} J e^{-i H_0 t}/dt=ie^{i H t} \big( H J - J H_0) e^{-i H_0 t},
$$
 it is desirable to find such $J$ that  the {\it effective perturbation} $HJ - J H_0$ be 
short-range (a pseudodifferential operator of order $-1-\varepsilon$ for $\varepsilon>0$).
This means that
  $\psi(x, \xi) = e^{i x \cdot \xi} p(x, \xi)$ is an approximate eigenfunction of the operator $H$,
with eigenvalue $|\xi|^2$, for each $\xi$. 
 Let us look for $\psi$ in the form $\psi(x, \xi) = e^{i\phi(x, \xi)}$. Then
\[
(-\Delta+ v(x)-|\xi|^2)\psi=(| \nabla \phi|^2 + v(x) - |\xi|^2 -i\Delta \phi)\psi,
\]
which leads to the eikonal equation  
$$
 | \nabla \phi|^2 + v(x) = |\xi|^2.
$$
Suppose that
\begin{equation}
\phi(x, \xi) = x \cdot \xi + \Phi(x, \xi),
\label{phPh}\end{equation}
where
\begin{equation}
\partial_x^\alpha \partial_\xi ^\beta \Phi(x,\xi) = O(|x|^{1-\rho-|\alpha|}).
\label{deriv-decay}\end{equation}
Then
\begin{equation}
 | \nabla \phi|^2 + v(x) - |\xi|^2=2 \xi \cdot \nabla \Phi+ | \nabla \Phi|^2 + v(x)
\label{nabla}\end{equation}
and neglecting the nonlinear term we obtain the equation
\begin{equation}
 2 \xi \cdot \nabla_x \Phi + v(x) = 0.
\label{nabla1}\end{equation}
We need its two solutions given by the equalities
\begin{equation}
\Phi_\pm(x, \xi) = \pm \frac1{2} \int_0^\infty \Bigl(v(x \pm t\xi) - v(\pm t\xi)\Bigr) \, dt.
\label{phi-eqn}\end{equation} 
Clearly, these functions $\Phi_\pm$   satisfy   assumption
\eqref{deriv-decay} but only off a conic neighbourhood of $x=\mp\xi$. 

Thus, two problems arise. The
first is that, to obtain in \eqref{J} a symbol $p_\pm$ from the class ${\cal S}^0_{\rho, \delta}$,
 we need to remove by a cut-off function $\zeta_\pm(x, \xi)$ a
  small conic neighbourhood  of the set $\hat{x} = \mp\hat{ \xi}$. Thus, we are obliged to
consider the wave operators for two different operators $J_\pm$. This idea appeared in \cite{IK3}
and will be realized below. The second problem is that the term $| \nabla \Phi|^2$ neglected in
\eqref{nabla} is `short range' only for $\rho>1/2$. Of course, it is easy to solve the eikonal equation by
iterations, considering at each step a linear equation of type \eqref{nabla1}, and to obtain a
`short range' error for arbitrary $\rho>0$. However, even after the cut-off by the function
$\zeta_\pm(x, \xi)$, we obtain the symbol $p_\pm$ from the class ${\cal S}^0_{\rho, 1-\rho}$, which
is bad if $\rho\leq 1/2$. To overcome this difficulty, we need to take into account the oscillating
nature of $p_\pm$ (see \cite{YP}).

\bigskip

{\bf 4}.
Below we suppose that $\rho>1/2$. Let $\sigma_\pm\in C^\infty$ be such that $\sigma_\pm(\theta)=1$
near $\pm 1$ and
$\sigma_\pm(\theta)=0$ near $\mp 1$.
We construct $J_\pm$    by the formula  \eqref{J} where
\begin{equation} 
p_\pm (x,\xi)=e^{ i\Phi_\pm (x,\xi) } 
\zeta_\pm (x,\xi)
\label{W1sy} \end{equation}
and the cut-off function $\zeta_\pm (x,\xi)$ essentially coincides with
 $\sigma_\pm ( \langle\hat{x},\hat{\xi}\rangle)$. We deliberately ignore here some technical
details which can be found in \cite{Y10}. For example, strictly speaking,   additional cut-offs of
low and high energies by a function of $ |\xi|^2 $ and   of a neighbourhood of $x=0$ by a function
of $x$ should be added to $\zeta_\pm (x,\xi)$.
  Then  the operators $J_\pm$ so constructed are pseudodifferential operators of
order $0$ and type $(\rho, \delta = 1 - \rho)$. 

Note that Theorem~\ref{smooth} extends automatically to the wave operators \eqref{WOJ}.
Thus, we are looking for a factorization
$$
H J_\pm - J_\pm H_0 = G^* \Omega G,
$$
where $G$ is locally $H_0$-  and   $H $-smooth and the operator $\Omega$ is bounded.   Let us recall  that the
operator $\ang{x}^{-r}$ is $H$-smooth, for any $r > 1/2$.  Since zeroth order pseudodifferential
operators are bounded on $L^2({\Bbb R}^d)$,  a factorization as above would have been true if $HJ -
JH_0$ were of order $-1 - \epsilon$ for some $\epsilon > 0$. 

However, the pseudodifferential operator $H J_\pm - J_\pm H_0$
has symbol decaying only as $|x|^{-1}$; 
 this is from one derivative of $-\Delta$ hitting $e^{ix \cdot \xi}$ and one hitting $\zeta_\pm$. 
Notice that the symbol of the `bad' term equals
\[
-2ie^{ i\Phi_\pm (x,\xi) } \langle \xi,\nabla_x \sigma_\pm  (
\langle\hat{x},\hat{\xi}\rangle)\rangle.
\]
It only decays as $|x|^{-1}$ but
 is supported outside a conic neighbourhood of the
set where $\hat{x} =  \hat{ \xi}$ or $\hat{x} = -\hat{ \xi}$.   Therefore,
$$
HJ - JH_0 = \sum_{j=1}^d (\ang{x}^{-1/2}\nabla^\perp_j)^* \Omega_j (\ang{x}^{-1/2}\nabla^\perp_j) + 
\ang{x}^{-r}\Omega_0 \ang{x}^{-r},\quad r>1/2,
$$
where $\Omega_j$, $j=0,1,\ldots,d$, are order zero pseudodifferential operators.   Using
Theorem~\ref{G} for the first $d$ terms, and the limiting absorption principle for the last one,
we obtain

\begin{prop}\label{Ex} 
 The wave operators $W_\pm(H, H_0, J_\tau)$ and $W_\pm(H_0, H , J_\tau^\ast)$ exist
  for $\tau=`+$' and $\tau=`-$'. 
\end{prop}

 The next step is to verify

\begin{prop}\label{Isom} 
 The operators $W_\pm(H,H_0;J_\pm )$ are isometric   and
$ W_\pm(H,H_0;J_\mp )=0$. 
\end{prop}

Indeed, it suffices to check that
\begin{equation}
  \slim_{t\rightarrow\pm\infty} (  J_\pm^\ast J_\pm -I) e^{-iH_0t}=0 
\label{WC1}\end{equation}
 and
\begin{equation}
  \slim_{t\rightarrow\pm\infty}  J_\mp ^\ast  J_\mp  e^{-iH_0t} =0.
\label{WC2}\end{equation}
According to  \eqref{W1sy},
 up to a compact term, $ J_\mp^\ast  J_\mp$ equals the pseudodifferential operator $Q_\mp$ 
with symbol $\zeta_\mp^2 (x,\xi) $.  If   $t\rightarrow\pm\infty$, then   the
stationary point $\xi=x/(2t)$ of the integral
\begin{equation}
  (Q_\mp e^{-iH_0t}f)(x)=(2\pi)^{-d/2}   \int_{{\Bbb R}^d}e^{i\langle \xi,x\rangle  -
i|\xi|^2t} \zeta_\mp^2 (x,\xi)
 \hat{f}(\xi) d\xi 
\label{WC3}\end{equation}
  does not belong to the support of the function
$\zeta_\mp^2 (x,\xi)$. Therefore supposing that  $f\in{\cal S}({\Bbb R}^d)$ and integrating by
parts, we estimate integral \eqref{WC3} by $C_N (1+|x|+|t|) ^{-N}$ for an arbitrary $N$. This
proves  \eqref{WC2}. To check \eqref{WC1}, we apply the same arguments to the PDO with symbol
$ \zeta_\pm^2(x,\xi) -1  $. 

Now it is easy to prove the asymptotic completeness.

\begin{thm}\label{WC3ZZ}
 Suppose that condition  \eqref{LR}   is fulfilled for $\rho >1/2$ and all $\alpha$.
 Then the wave operators $W_\pm (H,H_0; J_\pm)$ exist,  are isometric   and  complete.
\end{thm}

Since $ W_\pm(H_0,H;J_\mp^*)=  W_\pm^*(H,H_0;J_\mp ) $,
it follows from Proposition~\ref{Isom}  that $W_\pm(H_0,H;  J_\mp^\ast )=0$. This implies that 
\begin{equation}
\lim_{t\rightarrow\pm\infty} ||  J_\mp^\ast e^{-iHt}f|| =0,\quad f\in  {\cal H}^{ac}.
\label{eq:WC7}\end{equation} 
Let us choose the functions $\sigma_\pm$  in such a way that
$\sigma_+^2(\theta) + \sigma_-^2(\theta)=1$. Then 
\[
 J_+ J_+^\ast+J_-  J_-^\ast=\Id+K
\]
 for a compact operator $K$ and hence
\[
  ||  J_+^\ast e^{-iHt}f||^2+ ||  J_-^\ast e^{-iHt}f||^2=||f||^2 +o(1)
\]
as $|t|\rightarrow\infty$.
Now it follows from  \eqref{eq:WC7}   that
\[
\lim_{t\rightarrow\pm\infty}  ||   J_\pm^\ast e^{-iHt}f||= ||f||.
\] 
 This is equivalent to isometricity of the wave operators
\[
W_\pm(H_0,H;  J_\pm^\ast )=W_\pm^\ast(H,H_0;  J_\pm ),
\]
 so that         $W_\pm (H,H_0; J_\pm)$ are  complete.

\medskip

We emphasize that Theorem~\ref{WC3ZZ} and, essentially, its proof remain valid for an arbitrary
$\rho>0$.

Now it easy to justify the asymptotics \eqref{LRX}.  By existence and completeness of the wave
operator $W_\pm(H, H_0, J_\pm)$, we have
$$
\lim_{t\rightarrow\pm\infty}
\| e^{-i t H} f - J_\pm e^{-i t H_0} f_0^\pm \|= 0, \quad f_0^\pm =
W_\pm(H_0, H, J_\pm^\ast) f.
$$
The critical point $\xi(x, t)$ of the integral
\begin{equation} 
(J_\pm e^{-i t H_0} f_0^\pm)(x) =(2\pi)^{-d/2} \int_{{\Bbb R}^d} e^{i x \cdot \xi} e^{i \Phi_\pm
(x,\xi)} \zeta_\pm  (x, \xi) e^{-i t |\xi|^2} \hat f_0^\pm(\xi) \, d\xi 
\label{TD}\end{equation}
is defined by the equation
$$
2t\xi(x, t) = x + \nabla_\xi \Phi_\pm (x,\xi(x, t)),\quad \pm t>0,
$$
so that $\xi(x, t)=x/(2t)+ O(|t|^{  - \rho})$.
Applying stationary phase to the integral \eqref{TD}, we obtain formula \eqref{LRX}   with
function 
\[
\Xi(x,t)=  |x|^2/(4t) +\Phi_\pm (x,x/(2t)),
\]
which equals \eqref{WO6}.

\bigskip

{\bf 5}.
   In the long range case the scattering operator ${\bold S}$  is defined again by formula
\eqref{SOP} where
$W_\pm= W_\pm(H, H_0, J_\pm).$ 
Thus, again ${\bold S}$ is unitary and commutes with $H_0$, so defines a family of scattering
matrices
$S(\lambda)$ which are unitary operators on $L^2({\Bbb S}^{d-1})$ for $\lambda > 0$.

 However the structure of the spectrum of the scattering matrix in the short  and long range cases
are completely different. Typically, in the   long range case the spectrum of
$S(\lambda)$ covers the unit circle.  The nature of this spectrum is in general not known, except
in the radially symmetric case (that is, when the potential $v$ is a function only of $|x|$). In
that case, the scattering matrix commutes with the Laplacian on ${\Bbb S}^{d-1}$, so it breaks up
into the orthogonal sum of finite dimensional operators. Hence the spectrum is dense pure
point in this case. 
 
The kernel $s(\omega, \omega', \lambda)$ of $S(\lambda)$
 is still smooth away from the diagonal, but near the diagonal it, typically, has the form (for
$\rho<1$)
\begin{equation}
s(\omega, \omega', \lambda) \sim c(\omega,\lambda)|\omega - \omega'|^{-  (d-1)(1 + \rho^{-1})/2}
e^{i
\theta(\omega, \omega',\lambda)},
\label{Diag}\end{equation}
where  $\theta(\omega, \omega',\lambda)$ is
asymptotically homogeneous of order $1 - \rho^{-1}$ as $\omega - \omega'\rightarrow 0$. 

In the long range case practically nothing is known about eigenfunctions, which behave
asymptotically as plane waves. Moreover, as shows the explicit formula (see, e.g., \cite{LL}) for
the Coulomb potential $v(x)=v_0 |x|^{-1}$, in this case the separation of the asymptotics of
eigenfunctions into the sum of plane and spherical waves loses, to a large extent, its sense.
  On the contrary,  Theorem~\ref{aver} extends \cite{GY} to arbitrary long range potentials.

\section{The scattering matrix. High  energy and smoothness asymptotics}

This section relies on the paper \cite{Y1}.

\medskip

{\bf 1}.
Let us begin with the short range case, \eqref{power-decay} with $\rho > 1$. 
Then there is a stationary representation \eqref{sc-stat}  of the scattering matrix. Using
the Sobolev trace theorem and the dilation transformation $x\mapsto \lambda^{-1/2}x$, we can
show (see \cite{Y2}, for details) that, for any $r>1/2$, the operator \eqref{Gamma} satisfies
the estimate  
$$
\| \Gamma_0(\lambda) \ang{x}^{-r} \| \leq C(r) \lambda^{-1/4}.
$$
Similarly, one can control the dependence on $\lambda$ in the limiting absorption principle,
which yields
$$
\| \ang{x}^{-r} R(\lambda + i0) \ang{x}^{-r} \| \leq C(r) \lambda^{-1/2}.
$$
 
 The representation \eqref{sc-stat} allows us to study by a simple iterative procedure the
behaviour   of the scattering matrix  $S(\lambda)$ in
the two asymptotic regimes of interest, namely for high energies and in smoothness of its kernel.
In fact these two regimes are closely related.
 Namely, we use the resolvent identity \eqref{Res}
and substitute its right hand side in place of the resolvent $R(\lambda + i0)$ in \eqref{sc-stat}.
This gives us a series known as the Born approximation:
\[
S(\lambda) = \Id - 2\pi i \sum_{n = 0}^N (-1)^n \Gamma_0(\lambda) V  
\Big( R_0(\lambda + i0) V \Big)^n \Gamma_0^*(\lambda)  + \sigma_N(\lambda).
\]
The error term, $\sigma_N(\lambda)$ is $O(\lambda^{-(N-2)/2})$ in operator norm,
 and in addition gets smoother and smoother in the sense that $\sigma_N(\lambda) \in
\Sch_{\alpha(N)}$, where
$\alpha_N \to 0$ as $N \to \infty$.

However, the Born approximation has several drawbacks. First, it is actually quite complicated,
 involving multiple oscillating integrals of higher and higher dimensions as $n$ increases. Second,
it does not apply to long range potentials, or to magnetic, even short range, potentials. Here
we discuss another, much simpler,  form of approximation which applies
to all these situations but requires assumptions of the type \eqref{LR}. 

\bigskip

{\bf 2}.
To give an idea of the approach,
suppose first for simplicity that $v \in C_0^\infty(\RR^d)$ 
(though the argument can be applied to a wider class of potentials). Let the solution   of the
Schr\"odinger equation \eqref{efn-eqn} be defined by formula   
\eqref{psi}.
 By \eqref{sc-stat}, the integral kernel $k(\omega,\omega^\prime,\lambda)$ of the operator
$S(\lambda) -
\Id$ may be written as 
\begin{equation}
k(\omega,\omega^\prime ,\lambda)= -i \pi (2\pi)^{- d} \lambda^{(d-2)/2} \int_{{\Bbb R}^d} e^{-i
\sqrt{\lambda} x
\cdot\omega} v(x) \psi(x,   \omega', \lambda) \, dx.
\label{Kk}\end{equation}
Since this is an integral over a compact region, to analyze the
asymptotics of $k(\omega,\omega^\prime,\lambda)$ as $\lambda \to \infty$, it suffices to construct
the asymptotics of $\psi(x,   \omega,\lambda)$ for bounded $x$ as $\lambda \to \infty$.

 This can be
done by the following well known procedure \cite{Bu}. Assuming for a moment only
\eqref{LR} for $\rho>1$, one seeks
$\psi(x, \omega,\lambda)$ in the form
\begin{equation}
\psi(x, \xi) = e^{i x \cdot \xi}\text{b}(x, \xi),\quad \text{b}(x, \xi)=\sum_{n=0}^N (2i
|\xi|)^{-n} b_n(x, \hat\xi),
\quad b_0 = 1,
\quad \xi = \sqrt{\lambda} \omega.
\label{psi-ansatz}\end{equation}
Plugging \eqref{psi-ansatz} 
into the Schr\"odinger equation \eqref{efn-eqn}, and equating powers of $|\xi|$,
we obtain equations
\begin{equation}
\hat \xi \cdot \nabla_x b_{n+1} = -\Delta b_n + v b_n.
\label{b-eqn}\end{equation}
The remainder term 
$$
r_N (x,\xi) = e^{-i x \cdot \xi} (-\Delta+v(x) - |\xi|^2)\psi(x, \xi)
$$
is then given by
$$
r_N (x,\xi)= (2i |\xi|)^{-N} \big( -\Delta b_N (x,\hat\xi) + v (x) b_N (x,\hat\xi) \big).
$$
Equations \eqref{b-eqn} can be explicitly solved:
$$
b_{n+1}(x, \hat \xi) = \int_{-\infty}^0 \Big( -\Delta b_n(x + t\hat \xi, \hat \xi)
 + v(x + t\hat \xi) b_n(x + t\hat \xi, \hat \xi) \Big) dt.
$$
It is easy to see that
\begin{equation}\begin{aligned}
\big| \partial_x^\alpha \partial_\xi^\beta b_n(x, \hat \xi) \big| &\leq C_{\alpha, \beta} 
\ang{x}^{-(\rho - 1)n - |\alpha|} |\xi|^{-|\beta|}, \\
\big| \partial_x^\alpha \partial_\xi^\beta r_N(x,   \xi) \big| &\leq C_{\alpha, \beta}
\ang{x}^{-1-(\rho - 1)(N+1) - |\alpha|} |\xi|^{-N-|\beta|},
\end{aligned}\label{br}\end{equation}
except on arbitrary conic neighbourhoods of the bad direction $\hat x = \hat \xi$.

 If $v$ is compactly supported, then the estimates in $x$ are inessential, so it follows
from \eqref{Kk} and \eqref{psi-ansatz} that 
$$\begin{gathered}
k(\omega, \omega',\lambda) = -i\pi (2\pi)^{-  d } \lambda^{(d-2)/2} \sum_{n=0}^N
(2i \sqrt{\lambda})^{-n} \int_{{\Bbb R}^d} e^{i
\sqrt{\lambda} x \cdot (\omega^\prime-\omega )} v(x) b_n(x, \omega') \, dx \\
 + O(\lambda^{(d-3)/2 -N/2}).
\end{gathered}$$
Since $N$ is arbitrary, this gives the asymptotic expansion of the scattering amplitude as
$\lambda\rightarrow\infty$. We note that $k \in C^\infty({\Bbb S}^{d-1}\times
{\Bbb S}^{d-1})$  in the variables $\omega$ and $\omega^\prime$, so the smoothness asymptotics in
the case
$v\in C_0^\infty({\Bbb R}^d)$ is trivial.

\bigskip

{\bf 3}.
Finally, we give a universal formula which applies in both the long range and magnetic
cases. We shall need two approximate solutions $\psi_\pm$ of the Schr\"odinger equation 
$$
  \big( i \nabla + a(x) \big)^2\psi_\pm + v(x)\psi_\pm=|\xi|^2\psi_\pm.
$$
We suppose that a vector (or magnetic) potential $a(x)$ as well as scalar potential $v(x)$
satisfy the condition \eqref{LR} for some $\rho>0$. 
Let us look for $\psi_\pm$ in the form
\begin{equation}
\psi_\pm(x, \xi) = e^{i\phi_\pm(x, \xi)} \text{b}_\pm(x, \xi),
\label{psipm}\end{equation}
where $\phi=\phi_\pm$ is defined by formula \eqref{phPh} and $\Phi=\Phi_\pm$ satisfies
\eqref{deriv-decay}. 
Plugging 
\eqref{psipm} into the Schr\"odinger equation, we obtain the eikonal equation for $\phi$:
\begin{equation}
 |\nabla_x \phi|^2 - 2 a(x) \cdot \nabla_x \phi + v_0(x) = |\xi|^2, \quad
v_0(x) =   |a(x)|^2 + v(x). 
\label{eikonal}\end{equation}
If $\rho   > 1$ and $a = 0$, then one can set
$\Phi = 0$. In this case $\psi_-(x, \xi)=\psi(x, \xi)$ and $\psi_+(x, \xi)=\overline{\psi(x,
-\xi)}$ where the function $\psi(x, \xi)$ was constructed in the previous subsection.
 However, even if
$a$ is short range (but does not vanish), then, for the study of the limit
$\lambda\rightarrow\infty$, one cannot avoid the eikonal equation.

Once again, the equation \eqref{eikonal} for function $\phi=\phi_\pm$ defined by formula
\eqref{phPh} where $\Phi=\Phi_\pm$ can be solved by successive
approximations:
\[ 
\Phi (x, \xi) = \sum_{n=0}^{N_0} (2|\xi|)^{-n} \phi_n(x, \hat{\xi}).
\]
 Here   
\begin{equation*}\begin{gathered}
\hat \xi \cdot \nabla \phi_0 + \hat \xi \cdot a = 0,
 \\
\hat \xi \cdot \nabla \phi_1 + |\nabla \phi_0|^2 - 2 a \cdot \nabla \phi_0 + v_0 = 0
 \\
\hat \xi \cdot \nabla \phi_{n+1} + \sum_{m=0}^n \nabla \phi_m \cdot \nabla \phi_{n-m} -
2 a
\cdot \nabla \phi_n = 0 , \quad n \geq 1.
\end{gathered}\end{equation*}
So at every step we have an equation
$$
\hat \xi \cdot \nabla_x \phi(x,\hat{\xi}) + f(x,\hat{\xi}) = 0,
$$
with, possibly a long range function $f$. This equation can be solved by one of the two formulas
(cf. \eqref{nabla1}, \eqref{phi-eqn})
\begin{equation}
\phi_\pm(x, \xi) = \pm \int_0^\infty \Bigl(f(x \pm t \hat \xi,\hat{\xi}) - f(\pm t \hat
\xi,\hat{\xi})\Bigr)\, dt.
\label{eik-soln}\end{equation}
Using both signs `$+$' and `$-$',
 we obtain functions $\phi_\pm$ satisfying \eqref{eikonal} up to
a term $q_\pm(x, \xi)$ such that  
\[
\big| \partial_x^\alpha \partial_\xi^\beta q_\pm(x,   \xi) \big|
 \leq C_{\alpha, \beta} \ang{x}^{- N_0\rho - |\alpha|} |\xi|^{-N_0-|\beta|} 
\]
for all $(x, \xi)$ excluding an arbitrary conic neighbourhood of the direction 
$\hat x  =  \hat \xi$ for the minus sign, or $\hat x  = - \hat \xi$ for the plus sign.
One chooses and fixes $N_0$ such that $  N_0 \rho > 1$. 

Then for $\text{b}_\pm$ one has the transport
equation
\begin{equation}
  -2i   \xi\cdot\nabla \text{b}_\pm  +2i ( a-\nabla\Phi_\pm)\cdot\nabla \text{b}_\pm  
-\Delta\text{b}_\pm + ( -i \Delta \Phi_\pm  
 + i\text{div}\, a    +q_\pm ) \text{b}_\pm= 0.
\label{transport}\end{equation}
As before, one looks for $\text{b}_\pm=\text{b}_\pm^{(N)}$ in the form \eqref{psi-ansatz} which
gives standard equations 
\[
\hat\xi \cdot \nabla_x b_{n+1}^{(\pm)}(x,\hat\xi)= f_n^{(\pm)}(x,\hat\xi),
\]
where $f_n^{(\pm)}$ are determined by functions $b_1^{(\pm)},\ldots, b_n^{(\pm)}$.
All these equations are solved by formula \eqref{eik-soln} (but the term $f_n^{(\pm)}(\pm t\hat \xi)$ can be
dropped). Then the functions $b_n^{(\pm)}$ and the remainder $r_N^{(\pm)}$ in the transport
equation satisfy estimates of the type \eqref{br} for some $\rho>1$. As
$N\rightarrow\infty$, we obtain a function \eqref{psipm} satisfying the Schr\"odinger equation
with an arbitrary given accuracy both in the variables $x$ and $\xi$.

Now we can give an approximate formula for the scattering amplitude.
Remark that away from the diagonal,  it is not hard to show that the kernel $s(\omega, \omega',
\lambda)$ of the scattering matrix $S(\lambda)$ is smooth and $O(\lambda^{-\infty})$, so the main
point is to understand the kernel when
$\omega$ and $\omega'$ are close to some given point $\omega_0$.   Let
$\omega_0
\in
{\Bbb S}^{d-1}$ be arbitrary,  let $\Pi_{\omega_0}$ be the plane orthogonal to $\omega_0$ and
$x = z \omega_0 + y$, $y \in \Pi_{\omega_0}$.  
Set
\begin{equation}\begin{gathered}
s_0(\omega, \omega', \lambda) = \mp \pi i \lambda^{(d-2)/2}
 (2\pi)^{-d} \int_{\Pi_{\omega_0}} \Bigg( \overline{ \psi_+(y, \sql \omega) } \partial_z \psi_-(y,
\sql
\omega') \\ - \overline{ \partial_z \psi_+(y, \sql \omega) } \psi_-(y, \sql \omega') - 2i (a(y)
\cdot \omega_0) \overline{ \psi_+(y, \sql \omega) } \psi_-(y, \sql \omega') \Bigg) \, dy.
\end{gathered}\label{S_0}\end{equation}
Here $s_0 = s_0^{(N)}$   depends on $N$, since the construction above depends on a choice of $N$.
On the contrary, it does not depend on $N_0$ which is fixed. Let 
$\Omega_\pm=\Omega_\pm (\omega_0,\delta)\subset {\Bbb S}^{d-1}$ be determined by the condition
$\pm \omega \cdot
\omega_0 > \delta > 0$ and $\Omega=\Omega_+\cup\Omega_-$.
Then we have 

\begin{thm}\label{scatt-sing}
 For all integer  $p $, there is an $N=N(p)$ such that 
$$ 
\tilde s_0^{(N)}(\omega, \omega', \lambda):=s(\omega, \omega', \lambda) - s_0^{(N)}(\omega,
\omega',
\lambda)  \in C^p(\Omega \times \Omega) 
$$
and
$$
\big\| \tilde s_0^{(N)}(\cdot, \cdot, \lambda) \big\|_{C^p} \leq C \lambda^{-p}.
$$
\end{thm}

We make some brief comments on the proof of this theorem. Note that
 the representation formula \eqref{sc-stat}   holds   in the general short range case,
when both electric and magnetic potentials are present. However it is useless for the proof of
Theorem~\ref{scatt-sing} even for purely electric short range potentials.

As in the previous lecture, one considers instead
modified wave operators
$W_\pm(H,\\ H_0, J_\pm)$,
 where
$$
(J_\pm f)(x) = (2\pi)^{-d/2} \int_{{\Bbb R}^d} \psi_\pm(x, \xi) \zeta_\pm(x, \xi) \hat f(\xi) \,
d\xi
$$
and the functions $\psi_\pm(x, \xi)$ are defined by formula \eqref{psipm}.
Compared to the previous lecture, there are two important differences in the construction of the
operators $J_\pm$. First, we cannot neglect high energies and therefore have to control the
dependence on $\xi$ in all estimates. Second, for the proof of the existence and completeness of
the wave operators, it was sufficient to take $\text{b}_\pm=1$ in \eqref{psipm}. On the contrary,
now
$\text{b}_\pm= \text{b}_\pm^{(N)}$ is an approximate solution of the transport equation
\eqref{transport} depending on the parameter $N$,  which allows us to obtain an arbitrary good
approximation
$\psi_\pm=\psi_\pm^{(N)}$ to solutions of the Schr\"odinger equation.

 Writing $T_\pm$ for the effective perturbation,
$$
T_\pm = H J_\pm - J_\pm H_0,
$$
we have (see \cite{IK, Y10, Y2})
$$
S(\lambda) = S_0(\lambda) + S_1(\lambda),
$$
where
\begin{equation}\begin{aligned}
S_0(\lambda) &= -2 \pi i \Gamma_0(\lambda) J_+^* T_- \Gamma_0^*(\lambda),
 \\
S_1(\lambda) &=  2 \pi i \Gamma_0(\lambda) T_+^* R(\lambda + i0)  T_-  \Gamma_0^*(\lambda). 
\end{aligned}\label{S01}\end{equation}
Note that $T_\pm$ is a pseudodifferential operator of order $-1$,
 so that precise meaning of \eqref{S01} needs to be explained. The special reason why $S_0$ is
correctly defined is that the amplitude of the pseudodifferential operator
 $J_+^* T_-$ is zero in a neighbourhood of the set
where $\hat{x}$ is close to $\hat{\xi}$ or $-\hat{\xi}$  or, to put it differently, in a
neighbourhood of the conormal bundle to every sphere $  |\xi| =\sql  $. 

The term $S_1(\lambda)$ turns out to be negligible for large $N$. 

\begin{thm}
 For all integer  $p $, there is an $N=N(p)$ such that 
$$   
s_1^{(N)}(\omega, \omega', \lambda)   \in C^p({\Bbb S}^{d-1} \times {\Bbb S}^{d-1}),
$$
 and
$$
\big\| s_1^{(N)}(\cdot, \cdot, \lambda) \big\|_{C^p} \leq C \lambda^{-p}.
$$
\end{thm}

The proof relies on propagation estimates (see \cite{Mo2,JMP,J}) following from the Mourre
estimate  \eqref{Mourre}.
 We give an example of such an estimate. Let again $A$ be the generator of dilations. Then
for all integers $p$ 
\[
\big\| \ang{x}^p E_A({\Bbb R}_-) R(\lambda + i0) E_A({\Bbb R}_+) \ang{x}^p \big\| =
O(\lambda^{-1}). 
\]

Thus, the interesting part of the scattering matrix is contained in the term $S_0(\lambda)$
 of \eqref{S01}. It is very explicit representation, but has a drawback because it depends on the
cutoffs $\zeta_\pm$. So one has to transform the expression for $S_0(\lambda)$ to the invariant
expression \eqref{S_0}, which does not contain the cutoffs $\zeta_\pm$. This is the least obvious
part of the proof of Theorem~\ref{scatt-sing}. 

Finally, we note that  formula \eqref{Diag} for the diagonal singularity of the
scattering amplitude can be obtained applying the stationary phase method to integral \eqref{S_0}.

\end{document}